\documentclass[english,10pt, conference, letter]{IEEEtran}
\usepackage[T1]{fontenc}
\usepackage[latin9]{inputenc}
\usepackage{amsmath}
\usepackage{amsthm}
\usepackage{amssymb}
\usepackage{graphicx}

\makeatletter
\theoremstyle{plain}
\newtheorem{thm}{\protect\theoremname}
\theoremstyle{plain}
\newtheorem{prop}[thm]{\protect\propositionname}
\theoremstyle{plain}
\newtheorem{lem}[thm]{\protect\lemmaname}

\IEEEoverridecommandlockouts 
\newtheorem{assume}{Assumption}
\newtheorem{algor}{Algorithm}
\usepackage{cite}
\usepackage{hyperref}
\usepackage[margin=0.75in]{geometry}

\makeatother

\usepackage{babel}
\providecommand{\lemmaname}{Lemma}
\providecommand{\propositionname}{Proposition}
\providecommand{\theoremname}{Theorem}

\begin{document}

\title{\vspace{0.25in}
\LARGE \bf Modified Interior-Point Method for Large-and-Sparse Low-Rank
Semidefinite Programs}

\author{Richard Y. Zhang and Javad Lavaei
\thanks{This work was supported by the ONR YIP Award, DARPA YFA Award, AFOSR YIP Award, NSF CAREER Award, and NSF EPCN Award.}
\thanks{R.Y. Zhang and J. Lavaei are with the Department of Industrial Engineering and Operations Research, University of California, Berkeley, CA 94720, USA    {\tt\small ryz@berkeley.edu} and {\tt\small lavaei@berkeley.edu}}}
\maketitle
\begin{abstract}
Semidefinite programs (SDPs) are powerful theoretical tools that have
been studied for over two decades, but their practical use remains
limited due to computational difficulties in solving large-scale,
realistic-sized problems. In this paper, we describe a modified interior-point
method for the efficient solution of large-and-sparse low-rank SDPs,
which finds applications in graph theory, approximation theory, control
theory, sum-of-squares, etc. Given that the problem data is large-and-sparse,
conjugate gradients (CG) can be used to avoid forming, storing, and
factoring the large and fully-dense interior-point Hessian matrix,
but the resulting convergence rate is usually slow due to ill-conditioning.
Our central insight is that, for a rank-$k$, size-$n$ SDP, the Hessian
matrix is ill-conditioned only due to a rank-$nk$ perturbation, which
can be explicitly computed using a size-$n$ eigendecomposition. We
construct a preconditioner to ``correct'' the low-rank perturbation,
thereby allowing preconditioned CG to solve the Hessian equation in
a few tens of iterations. This modification is incorporated within
SeDuMi, and used to reduce the solution time and memory requirements
of large-scale matrix-completion problems by several orders of magnitude.
\end{abstract}

\section{Introduction}

\global\long\def\R{\mathbb{R}}
\global\long\def\C{\mathbb{C}}
\global\long\def\S{\mathbb{S}}
\global\long\def\A{\mathbf{A}}
\global\long\def\B{\mathbf{B}}
\global\long\def\H{\mathbf{H}}
\global\long\def\P{\mathbf{P}}
\global\long\def\D{\mathbf{D}}
\global\long\def\Sch{\mathbf{S}}
\global\long\def\Z{\mathbf{Z}}
\global\long\def\U{\mathbf{U}}
\global\long\def\Real{\mathrm{Re}\,}
\global\long\def\Imag{\mathrm{Im}\,}
\global\long\def\tr{\mathrm{tr}\,}
\global\long\def\cond{\mathrm{cond}\,}
\global\long\def\diag{\mathrm{diag}\,}
\global\long\def\rank{\mathrm{rank}\,}
\global\long\def\cone{\mathcal{K}}
\global\long\def\nei{\mathcal{N}_{-\infty}}
\global\long\def\vector{\mathrm{vec}\,}
\global\long\def\mat{\mathrm{mat}\,}
Consider the size-$n$ semidefinite program with $m$ constraints
\begin{align}
X^{\star}=\text{ minimize } & C\bullet X\tag{SDP}\label{eq:SDP}\\
\text{subject to } & A_{i}\bullet X=b_{i}\quad\forall i\in\{1,\ldots,m\}\nonumber \\
 & X\succeq0,\nonumber 
\end{align}
and its Lagrangian dual
\begin{align}
\{y^{\star},S^{\star}\}=\text{ maximize } & b^{T}y\tag{SDD}\label{eq:SDD}\\
\text{subject to } & \sum_{i=1}^{m}y_{i}A_{i}+S=C\nonumber \\
 & S\succeq0.\nonumber 
\end{align}
Each matrix is $n\times n$ real symmetric (an element of $\S^{n}$);
$\bullet$ denotes the associated matrix inner product $A\bullet B=\tr A^{T}B$;
and $X\succeq0$ and $S\succ0$ ($X\in\S_{+}^{n}$ and $S\in\S_{++}^{n}$)
indicate that $X$ is symmetric positive semidefinite and $S$ is
symmetric positive definite. In case of nonunique solutions, we use
$\{X^{\star},y^{\star},S^{\star}\}$ to refer to the \emph{analytic
center} of the solution set.

In this paper, we consider large-and-sparse low-rank SDPs, for which
the number of nonzeros in the data $A_{1},\ldots,A_{m}$ is small,
and $k\triangleq\rank X^{\star}$ is known \emph{a priori} to be very
small relative to the dimensions of the problem, i.e. $k\ll n$. Such
problems widely appear as the convex relaxations of ``hard'' optimization
problems in graph theory~\cite{sojoudi2014exactness}, approximation
theory~\cite{candes2010power,candes2012exact,lasserre2009moments},
control theory~\cite{valmorbida2016stability,lasserre2009moments,zhang2016thesis},
and power systems~\cite{madani2015admm,madani2015convex}. They are
also the fundamental building blocks for global optimization techniques
based upon polynomial sum-of-squares~\cite{parrilo2003semidefinite}
and the generalized problem of moments~\cite{lasserre2009moments}.

Interior-point methods are the most reliable approach for solving
small- and medium-scale SDPs, but become prohibitively time- and memory-intensive
for large-scale problems. A fundamental issue is their inability to
exploit problem structure, such as the sparsity of the data and the
low-rank feature of the solution, to substantially reduce complexity.
In other words, interior-point methods solve highly sparse, rank-one
SDPs in approximately the same time as dense, full-rank SDPs of the
same size. 

In this paper, we present a modification to the standard interior-point
method that makes it substantially more efficient for large-and-sparse
low-rank SDPs. More specifically, our algorithm solves a rank-$k$
SDP in $\Theta(n^{3}k^{3})$ time and $\Theta(n^{2}k^{2})$ memory,
under some mild nondegeneracy and sparsity assumptions. In Section~\ref{sec:Numerical-Results},
we give numerical results to show that our method is up to a factor
of $n$ faster than the standard interior-point method for problems
with $m\sim n$ constraints, and up to a factor of $n^{3}$ faster
for problems with $m\sim n^{2}$ constraints.

\subsection{Assumptions}

We begin with some nondegeneracy assumptions, which are standard for
interior-point methods.

\begin{assume}[Nondegeneracy]\label{ass:nondegen}We assume:
\begin{enumerate}
\item (Slater's condition) There exist $X\succ0$, $y$, and $S\succ0$,
such that $A_{i}\bullet X=b_{i}$ and $\sum_{i}y_{i}A_{i}+S=C$.
\item (Strict complementarity) $\rank(X^{\star})+\rank(S^{\star})=n$.
\end{enumerate}
\end{assume}

These are generic properties of SDPs, and are satisfied by almost
all instances~\cite{alizadeh1997complementarity}. Note that Slater's
condition is satisfied in solvers like SeDuMi~\cite{sturm1999sedumi}
and MOSEK~\cite{mosek2015} using the homogenous self-dual embedding
technique~\cite{ye1994homogeneous}.

We further assume that the data matrices $A_{1},\ldots,A_{m}$ are
structured in a way that allow certain matrix-implicit operations
to be efficiently performed. 

\begin{assume}[Sparsity]\label{ass:sparse}Define the matrix $\A\triangleq[\vector A_{1},\ldots,\vector A_{m}]$.
We assume that matrix-vector products with $\A$, $\A^{T}$ and $(\A^{T}\A)^{-1}$
may each be applied in $O(m)$ flops and memory.

\end{assume}

Versions of this assumption appear in most large-scale SDP algorithms,
spanning both first-order~\cite{wen2010alternating,madani2015admm,kalbat2015fast,o2016conic,zhang2016convergence,zhang2016thesis}
and second-order methods~\cite{toh2002solving,zhao2010newton}. The
assumption is satisfied by any sparse data whose normal matrix $\A^{T}\A$
admits a sparse Cholesky factorization. 

\subsection{Related work}

The desire to effectively exploit problem structure in large-scale
SDPs has motivated a number of algorithms. It is convenient to categorize
them into three distinct groups: 

The first group is based on using sparsity in the data to decompose
the size-$n$ conic constraint $X\succeq0$ into many smaller conic
constraints over submatrices of $X$. In particular, when the matrices
$C,A_{1},\ldots,A_{m}$ share a common sparsity structure with a chordal
graph with bounded treewidth $\tau$, a technique known as \emph{chordal
decomposition} or \emph{chordal conversion} can be used to reformulate
(\ref{eq:SDP})-(\ref{eq:SDD}) into a problem containing only size-$(\tau+1)$
semidefinite constraints~\cite{fukuda2001exploiting}; see also~\cite{vandenberghe2015chordal}.
While the technique is only applicable to chordal SDPs with bounded
treewidths, it is able to reduce the cost of a size-$n$ SDP all the
way down to the cost of a size-$n$ linear program, sometimes as low
as $O(\tau^{3}n)$. Indeed, chordal sparsity can be guaranteed in
many important applications~\cite{vandenberghe2015chordal,madani2015convex},
and software exist to automate the chordal reformulation~\cite{kim2011exploiting}.

The second group is based on applying first-order methods for nonlinear
programming, such as conjugate gradients~\cite{toh2002solving,zhao2010newton,zhang2016thesis}
and ADMM~\cite{wen2010alternating,madani2015admm,kalbat2015fast,o2016conic,zhang2016convergence},
either to (\ref{eq:SDP}) directly, or to the Newton subproblem associated
with an interior-point solution of (\ref{eq:SDP}). These algorithms
have inexpensive per-iteration costs but a sublinear worst-case convergence
rate, computing an $\epsilon$-accurate solution in $O(1/\epsilon)$
time. They are most commonly used to solve very large-scale SDPs to
modest accuracy.

The third group is based on the outer product factorization $X=RR^{T}$.
These methods use the low-rank of $X^{\star}$ to reduce the number
of decision variables in (\ref{eq:SDP}) from $\sim n^{2}$ to $nk$~\cite{burer2003nonlinear,journee2010low}.
The problem being solved is no longer convex, so only local convergence
can be guaranteed. Nevertheless, time and memory requirements are
substantially reduced, and these methods have been used to solve very
large-scale low-rank SDPs to excellent precision; see the computation
results in~\cite{burer2003nonlinear,journee2010low}.

Our method is similar in spirit to methods from the second group,
but makes much stronger convergence guarantees. More specifically,
we guarantee that the method converges globally to $\{X^{\star},y^{\star},S^{\star}\}$
at a linear rate, producing an $\epsilon$-accurate solution in $O(\log(1/\epsilon))$
time. At the same time, the method remains applicable for SDPs that
are sparse but nonchordal. Indeed, in Section~\ref{sec:Numerical-Results},
we present strong computational results for the matrix completion
problem, which cannot be efficiently solved using methods from the
first group. We mention, however, that the method has a higher memory
requirement than methods from the third group, due to its need to
explicitly store the matrix variables $X$ and $S$.

\subsection{Notations}

Most of our notations are standard except the following. Given a positive
definite matrix $X\in\S_{++}^{n}$, we order its eigenvalues $\lambda_{1}(X)\ge\cdots\ge\lambda_{n}(X)$,
and define its condition number $\kappa(X)=\lambda_{1}(X)/\lambda_{n}(X)$.
We sometimes use $\lambda_{\max}(X)\equiv\lambda_{1}(X)$ and $\lambda_{\min}(X)\equiv\lambda_{n}(X)$
for emphasis. We use ``$\vector$'' and ``$\otimes$'' to refer
to the (nonsymmetricized) vectorization and Kronecker product, which
satisfy the identity $\vector AXB^{T}=(A\otimes B)\vector X$. We
use $\diag(A,B)=\left[\begin{smallmatrix}A & 0\\
0 & B
\end{smallmatrix}\right]$ to refer to the matrix direct sum.

\section{Interior-Point Methods}

Consider replacing the nonsmooth, convex constraint $X\succeq0$ in
(\ref{eq:SDP}) by the smooth, strongly convex, and self-concordant
penalty function $\mu\log\det X$, as in
\begin{align}
X_{\mu}=\text{ minimize } & C\bullet X-\mu\log\det X\tag{SDP\ensuremath{\mu}}\label{eq:PotPrim}\\
\text{subject to } & A_{i}\bullet X=b_{i}\;\forall i\in\{1,\ldots,m\}.\nonumber 
\end{align}
The resulting problem has Lagrangian dual
\begin{align}
\{y_{\mu},S_{\mu}\}=\text{ maximize } & b^{T}y+\mu\log\det S\tag{SDD\ensuremath{\mu}}\label{eq:PotDual}\\
\text{subject to } & \sum_{i=1}^{m}y_{i}A_{i}+S=C.\nonumber 
\end{align}
For different values of $\mu>0$, the corresponding solutions $\{X_{\mu},y_{\mu},S_{\mu}\}$
define a trajectory in the feasible region of (\ref{eq:SDP})-(\ref{eq:SDD})
that approaches $\{X^{\star},y^{\star},S^{\star}\}$ as $\mu\to0^{+}$.
This trajectory is known as the \emph{central path}, and $\mu$ is
known as the duality gap parameter, because $n\mu=X_{\mu}\bullet S_{\mu}=C\bullet X_{\mu}-b^{T}y_{\mu}$
is the duality gap of the feasible point $\{X_{\mu},y_{\mu},S_{\mu}\}$
in (\ref{eq:SDP})-(\ref{eq:SDD}).

All interior-point methods work by using Newton's method to approximately
solve (\ref{eq:PotPrim}), (\ref{eq:PotDual}), or their joint Karush\textendash Kuhn\textendash Tucker
(KKT) equations, while making decrements in the duality gap parameter
$\mu$. Most modern SDP solvers are of the path-following type, and
explicitly keep their iterates within a feasible neighborhood of the
central path
\begin{align}
\mathcal{N}_{\infty}^{-}(\gamma) & \triangleq\left\{ \{X,y,S\}\text{ feas.}:\lambda_{\min}(XS)\ge\frac{\gamma}{n}\tr XS\right\} ,\label{eq:nei}
\end{align}
where $\gamma\in(0,1)$ quantifies the ``size'' of the neighborhood.
The resulting interior-point method has a formal iteration complexity
of $O(n\log\epsilon^{-1})$, but always converges within tens of iterations
in practice~\cite[Ch.5]{wright1997primal}. 

Each iteration of an interior-point method solves 1-3 quadratic approximations
of (\ref{eq:PotDual})
\begin{align}
\text{ maximize } & b^{T}y-\frac{1}{2}\|W^{\frac{1}{2}}(S-Z)W^{\frac{1}{2}}\|_{F}^{2}\label{eq:Newt}\\
\text{subject to } & \sum_{i=1}^{m}y_{i}A_{i}+S=C,\nonumber 
\end{align}
in which $W,Z\in\S_{++}^{n}$ are used by the algorithm to approximate
the log-det penalty. Substituting $S=C-\sum_{i=1}^{m}y_{i}A_{i}$
into the objective (\ref{eq:Newt}) yields an unconstrained problem
with first-order optimality conditions:
\begin{equation}
A_{i}\bullet\left[W\left(\sum_{j=1}^{m}y_{j}A_{j}\right)W\right]=\underbrace{b_{i}+A_{i}\bullet W(C-Z)W}_{r_{i}}\label{eq:hesseqn2}
\end{equation}
for all $i\in\{1,\ldots,m\}$. Vectorizing the matrix variables allows
(\ref{eq:hesseqn2}) to be compactly written as
\begin{equation}
\underbrace{(\A^{T}\D\A)}_{\H}y=r\label{eq:hesseqn1}
\end{equation}
where $\A=[\vector A_{1},\ldots,\vector A_{m}]$ and $\D=W\otimes W$.
Once $y$ is computed, the variables $S=C-\sum_{i=1}^{m}y_{i}A_{i}$
and $X=W(Z-S)W$ are easily recovered.

Since the interior-point method converges in tens of iterations, the
cost of solving (\ref{eq:SDP})-(\ref{eq:SDD}) is essentially the
same as that of solving the Hessian equation $\H y=r$, up to a modest
multiplicative constant. Or put in another way, an interior-point
method can be thought of as a technique to convert the nonsmooth conic
problems (\ref{eq:SDP})-(\ref{eq:SDD}) into a small sequence of
unconstrained least-squares problems~\cite[Ch.11]{boyd2004convex}. 

\subsection{Solving the Hessian equation}

The computation bottleneck in every interior-point method is the solution
of the Hessian equation $\H y=r$. The standard approach found in
the vast majority of interior-point solvers is to form $\H$ explicitly
and to factor it using Cholesky factorization. An important feature
of interior-point methods for SDPs is that the matrix $W$ is fully-dense,
so the cost of forming and factoring the fully-dense $m\times m$
Hessian matrix $\H$ using dense Cholesky factorization is $O(n^{3}m+n^{2}m^{2}+m^{3})$
time and $\Theta(m^{2}+n^{2})$ memory.

Alternatively, the Hessian equation may be solved using an iterative
method like conjugate gradients (CG). We defer to standard texts~\cite{barrett1994templates}
for implementation details, and only note that the method requires
a single matrix-vector product with the governing coefficient matrix
at each iteration. In exact arithmetic, CG converges to the exact
solution of the Hessian equation $\H y=r$ within $m$ iterations,
thereby producing a complexity of $O(n^{3}m+n^{2}m^{2})$ time and
$\Theta(n^{2}+m)$ memory, which is strictly better than Cholesky
factorization. 

However, in finite precision, CG does not terminate in $m$ steps
due to the accumulation of round-off error. Instead, the method converges
linearly, with a convergence rate related to the condition number
of the governing matrix. 
\begin{prop}[{\cite[p.53]{greenbaum1997iterative}}]
\label{prop:CG}Given $x^{0},b\in\R^{n}$, $A\in\S_{++}^{n}$, define
$x^{\star}=A^{-1}b$. Then, the $i$-th iterate of CG generates satisfies
\begin{equation}
\frac{\|x^{i}-x^{\star}\|}{\|x^{0}-x^{\star}\|}\le2\sqrt{\kappa_{1}}\left(\frac{\sqrt{\kappa_{j}}-1}{\sqrt{\kappa_{j}}+1}\right)^{i-j}\label{eq:cgcheb}
\end{equation}
with condition numbers $\kappa_{j}=\lambda_{j}(A)/\lambda_{\min}(A)$
and $j\in\{1,\ldots,n\}$.
\end{prop}
The Hessian matrix $\H$ becomes increasing ill-conditioned as the
outer interior-point method makes progress towards the solution. Its
condition number scales $\kappa(\H)=O(1/\mu^{2})$, where $\mu$ is
the duality gap parameter at the current interior-point iteration.
This ill-conditioning gives any CG-based interior-point method a sublinear
worse-case time complexity, converging to an $\epsilon$-accurate
solution of (\ref{eq:SDP})-(\ref{eq:SDD}) in $O(1/\epsilon)$ time.

Instead, all successful CG-based solution of the Hessian equation
rely on an effective preconditioner, and a modification to CG named
\emph{preconditioned} conjugate gradients (PCG). Each PCG iteration
requires a single matrix-vector product with the governing matrix,
and a single \emph{solve} with the preconditioner; see e.g.~\cite{barrett1994templates,greenbaum1997iterative}.
\begin{prop}
\label{prop:PCG}Given $x^{0},b\in\R^{n}$ , $A\in\S_{++}^{n}$, and
preconditioner $P\in\S_{++}^{n}$, define $x^{\star}=A^{-1}b$. Then,
the $i$-th iterate of PCG generates satisfies (\ref{eq:cgcheb})
with $\kappa_{j}=\lambda_{j}(P^{-1}A)/\lambda_{n}(P^{-1}A)$.
\end{prop}
If a preconditioner $\tilde{\H}$ can be constructed to be spectrally
similar to $\H$ (in the specific sense described in Proposition~\ref{prop:PCG}),
then PCG allows us to solve a Hessian $\H y=r$ by solving a few instances
of the preconditioner equation $\tilde{\H}y=r$. 

\subsection{Ill-conditioning in the scaling matrix}

The matrix $W\in\S_{++}^{n}$ is known as the \emph{scaling matrix},
and captures the curvature of the log-det penalty function. Different
interior-point methods differ primarily how the scaling matrix $W$
is constructed. Given the current iterate $\{\hat{X},\hat{y},\hat{S}\}$,
we consider three types of scalings:
\begin{itemize}
\item \textbf{Primal scaling.} Set $W\gets\hat{X}$. Used in the original
projective conic interior-point method by Nesterov \& Nemirovski~\cite{nesterov1994interior}. 
\item \textbf{Dual scaling.} Set $W\gets\hat{S}^{-1}$. Used in the log-determinant
barrier method~\cite{vandenberghe1998determinant}. 
\item \textbf{Nesterov-Todd (NT) scaling.} Set $W$ to be the unique positive
definite matrix satisfying $\hat{X}=W\hat{S}W$. This is the most
widely used scaling method for semidefinite programming, found in
SeDuMi~\cite{sturm1999sedumi} and MOSEK~\cite{mosek2015}.
\end{itemize}
In all three cases, the scaling matrix $W$ becomes progressively
ill-conditioned as the interior-point method makes progress towards
the solution. This is the mechanism that causes the Hessian matrix
$\H$ to become ill-conditioned; see~\cite{alizadeh1997complementarity}.
\begin{lem}
\label{lem:XS_eig}Under Assumption~\ref{ass:nondegen}, fix $\mu_{0}>0$
and $\gamma\in(0,1)$. Then, for all points $\{X,y,S\}$ with
\[
\{X,y,S\}\in\mathcal{N}_{\infty}^{-}(\gamma),\qquad\mu\le\mu_{0}
\]
(where $\mu=\frac{1}{n}\tr XS$), there are constants $C_{0}$ and
$C_{1}$ such that
\begin{gather}
\lambda_{1}(X)\le C_{0},\qquad\lambda_{1}(S)\le C_{0},\label{eq:eig1}\\
\lambda_{k}(X)\ge C_{1}\gamma,\qquad\lambda_{n-k}(S)\ge C_{1}\gamma,\label{eq:eigk}\\
\lambda_{k+1}(X)\le\mu/C_{1},\qquad\lambda_{n-k+1}(S)\le\mu/C_{1},\label{eq:eigk1}\\
\lambda_{n}(X)\ge\gamma\mu/C_{0},\qquad\lambda_{n}(S)\ge\gamma\mu/C_{0}.\label{eq:eign}
\end{gather}
\end{lem}
\begin{IEEEproof}
This is the SDP version of Lemma 5.13 in \cite{wright1997primal},
which was stated for LPs. 
\end{IEEEproof}
\begin{prop}
\label{prop:Weig}Under the conditions in Lemma~\ref{lem:XS_eig},
let $W$ be the primal, dual, or NT scaling matrix computed from $\{X,S\}$.
Then,
\[
\frac{\lambda_{1}(W)}{\lambda_{k}(W)}=O(1),\;\frac{\lambda_{k+1}(W)}{\lambda_{n}(W)}=O(1),\;\frac{\lambda_{k}(W)}{\lambda_{k+1}(W)}=\Theta(1/\mu).
\]
\end{prop}
\begin{IEEEproof}
Lemma~\ref{lem:XS_eig} establishes these conditions for $X$ and
$S^{-1}$. For NT scaling, let us note that $W=X\#S^{-1}$, where
$\#$ is the (metric) geometric mean operator of Ando~\cite{ando1979concavity}.
Then, Ando's matrix arithmetic-geometric inequality implies $\frac{1}{2}(X+\mu S^{-1})\succeq X\#(\mu S^{-1})=\frac{1}{\sqrt{\mu}}W$
and $\frac{1}{2}(\mu X^{-1}+S)\succeq(\mu X^{-1})\#S=\sqrt{\mu}W^{-1}$. 
\end{IEEEproof}

\section{Preconditioning the Hessian Matrix}

In this section, we develop a preconditioner $\tilde{\H}$ that is
both easy to invert, and also serves as a good spectral approximation
for $\H$. More specifically, we prove that PCG with $\tilde{\H}$
as preconditioner solves the Hessian equation $\H y=r$ to machine
precision in a \emph{constant} number of iterations, irrespective
of $\mu$. 

\subsection{The main idea}

The preconditioner is based off the observation that the scaling matrix
$W$ becomes ill-conditioned only due to the presence of $k$ large
outlier eigenvalues. Using a single size-$n$ eigendecomposition,
$W$ can be decomposed into a well-conditioned component and a low-rank
perturbation, as in 
\begin{equation}
W=W_{0}+UU^{T},\label{eq:WoUUt}
\end{equation}
where $\kappa(W_{0})\in O(1)$ and $\rank U\le k$. Indeed, let us
partition the eigenvalues and eigenvectors of $W$ into two groups,
\begin{equation}
W=\begin{bmatrix}V_{s} & V_{\ell}\end{bmatrix}\begin{bmatrix}\Lambda_{s} & 0\\
0 & \Lambda_{\ell}
\end{bmatrix}\begin{bmatrix}V_{s} & V_{\ell}\end{bmatrix}^{T},\label{eq:Wsplit1}
\end{equation}
putting the smallest $n-k$ eigenvalues into $\Lambda_{s}$, and the
$k$ largest eigenvalues into $\Lambda_{\ell}$. Then, choosing any
$\tau$ to satisfy $\lambda_{\min}(\Lambda_{s})\le\tau<\lambda_{\max}(\Lambda_{s})$,
the following
\begin{equation}
W=\underbrace{\begin{bmatrix}V_{s} & V_{\ell}\end{bmatrix}\begin{bmatrix}\Lambda_{s} & 0\\
0 & \tau I
\end{bmatrix}\begin{bmatrix}V_{s} & V_{\ell}\end{bmatrix}^{T}}_{W_{0}}+\underbrace{V_{\ell}(\Lambda_{\ell}-\tau I)V_{\ell}^{T}}_{UU^{T}}\label{eq:Wsplit2}
\end{equation}
implements the desired splitting in (\ref{eq:WoUUt}).

Since $W_{0}$ is well-conditioned, it can be well approximated by
a scaled identity matrix. Substituting $W_{0}\approx\tau I$ in (\ref{eq:WoUUt})
yields a low-rank perturbation of the identity
\begin{equation}
\tilde{W}=\tau I+UU^{T}.\label{eq:tIUUt}
\end{equation}
Matrix-vector products with $\tilde{W}^{-1}$ can be efficiently performed
using the Sherman\textendash Morrison\textendash Woodbury formula
\begin{equation}
\tilde{W}^{-1}=(\tau I+UU^{T})^{-1}=\tau^{-1}I-\tau^{-1}US^{-1}U^{T},\label{eq:wsmf}
\end{equation}
in which $S=\tau I+U^{T}U$ is a $k\times k$ positive definite Schur
complement. By virtue of $\tau I$ being a good spectral approximation
of $W_{0}$, the matrix $\tilde{W}$ is also a good spectral approximation
for $W$.
\begin{lem}
\label{lem:kappaW}Let $W$ and $\tilde{W}$ be defined in (\ref{eq:WoUUt})
and (\ref{eq:tIUUt}), and choose $\lambda_{\min}(W_{0})\le\tau\le\lambda_{\max}(W_{0})$.
Then, $\kappa(W,\tilde{W})=\kappa(W_{0})$.
\end{lem}
\begin{IEEEproof}
Define $F\triangleq\begin{bmatrix}\sqrt{\tau}I_{n} & U\end{bmatrix}^{T}$,
so that $\tilde{W}=F^{T}F$ and $W=F^{T}\diag(\frac{1}{\tau}W_{0},I_{k})F$.
Define $Q\triangleq F(F^{T}F)^{-1/2}$, and observe that $Q$ is a
matrix with orthonormal columns. Then, $\tilde{W}^{-1/2}W\tilde{W}^{-1/2}=Q^{T}\diag(\frac{1}{\tau}W_{0},I_{k})Q$.
Applying the Cauchy interlacing eigenvalues theorem, we have $\kappa(W,\tilde{W})=\kappa(\tilde{W}^{-1/2}W\tilde{W}^{-1/2})\le\kappa\left(\diag(\frac{1}{\tau}W_{0},I_{k})\right)=\kappa(W_{0})$.
\end{IEEEproof}

\subsection{Extending to the Hessian matrix}

Similarly, the Hessian matrix $\H$ becomes ill-conditioned only due
to the presence of $nk$ large outlier eigenvalues. Substituting the
splitting (\ref{eq:WoUUt}) into $\H=\A^{T}\D\A$ yields
\begin{multline}
\H=\A^{T}(W_{0}\otimes W_{0}+UU^{T}\otimes W_{0}\\
+W_{0}\otimes UU^{T}+UU^{T}\otimes UU^{T})\A.\label{eq:hessExp1}
\end{multline}
The terms can be collected using the following observation.
\begin{lem}
\label{lem:sym_kron}For any $X,Y\in\R^{n\times n}$, not necessarily
symmetric, we have $\A^{T}(X\otimes Y)\A=\A^{T}(Y\otimes X)\A.$
\end{lem}
\begin{IEEEproof}
We have $[\A^{T}(X\otimes Y)\A]_{i,j}=\tr A_{i}XA_{j}Y^{T}=\tr A_{i}YA_{j}X^{T}=[\A^{T}(Y\otimes X)\A]_{i,j}$
due to the symmetry of $A_{i}$, $A_{j}$, and the cyclic property
of the trace operator.
\end{IEEEproof}
Applying Lemma~\ref{lem:sym_kron} yields a well-conditioned plus
low-rank splitting for the matrix $\H$, as in
\begin{equation}
\H=\underbrace{\A^{T}(W_{0}\otimes W_{0})\A}_{\H_{0}}+\underbrace{\A^{T}(U\otimes Z)(U\otimes Z)^{T}\A}_{\U\U^{T}}.\label{eq:hessExp2}
\end{equation}
where $Z$ is any matrix (not necessarily unique) satisfying $ZZ^{T}=2W_{0}+UU^{T}$. 

Again, we approximate the well-conditioned matrix $W_{0}$ using a
scaled identity. Substituting $W_{0}\approx\tau I$ yields
\begin{equation}
\tilde{\H}=\tau^{2}\A^{T}\A+\U\U^{T},\label{eq:AtIUUtAt}
\end{equation}
whose inverse can also be expressed using the Sherman\textendash Morrison\textendash Woodbury
formula
\begin{equation}
\tilde{\H}^{-1}=(\tau^{2}\A^{T}\A)^{-1}(I-\U\Sch^{-1}\U^{T}(\A^{T}\A)^{-1}),\label{eq:SMW2}
\end{equation}
with $\Sch=\tau^{2}I+\U^{T}(\A^{T}\A)^{-1}\U$. Note that each matrix-vector
product with $\U$ and its transpose can be efficiently performed
by exploiting the Kronecker structure,\begin{subequations}\label{eq:Umv}
\begin{align}
\U\vector X & =\A^{T}(U\otimes Z)\vector X=[\tr A_{i}(ZX)U^{T}]_{i=1}^{m},\\
\U^{T}y & =(U\otimes Z)^{T}\A y=Z^{T}\left(\sum_{i=1}^{m}y_{i}A_{i}\right)U,
\end{align}
\end{subequations}in $2n^{2}k$ flops and a call to $\A^{T}$ or
$\A$. Hence, (\ref{eq:SMW2}) can be efficiently evaluated assuming
that efficient matrix-vector products with $\A$, $\A^{T}$, and $(\A^{T}\A)^{-1}$
are available (Assumption~\ref{ass:sparse}). 

We can repeat the same arguments as before to show that $\tilde{\H}$
is a good spectral approximation of $\H$.
\begin{lem}
\label{lem:Htilde}Given $\H=\A^{T}\D\A$, let $\tilde{\H}$ be defined
in (\ref{eq:AtIUUtAt}). Choose $\tau$ to satisfy $\lambda_{\min}(W_{0})\le\tau\le\lambda_{\max}(W_{0})$.
Then, $\kappa(\H,\tilde{\H})\le\kappa^{2}(W_{0})$.
\end{lem}
\begin{IEEEproof}
Define $F\triangleq\begin{bmatrix}\tau\A^{T} & \U\end{bmatrix}^{T}$
and repeat the proof of Lemma~\ref{lem:kappaW}.
\end{IEEEproof}
In view of Lemma~\ref{lem:Htilde} and Proposition~\ref{prop:PCG},
we find that PCG with $\tilde{\H}$ as preconditioner solves the Hessian
equation $\H y=r$ in a constant number of iterations.

\subsection{Complexity analysis}

The full PCG solution procedure is summarized as Algorithm~\ref{alg:prec}.

\begin{algor}\label{alg:prec}Input: Right-hand side $r\in\R^{m}$,
relative accuracy $\epsilon>0$, scaling matrix $W\in\S_{++}^{n}$,
solution rank $k>0$, and efficient matrix-vector products with $\A$,
$\A^{T}$, and $(\A^{T}\A)^{-1}$.\\
Output: An $\epsilon$-accurate solution vector $y\in\R^{m}$ for
the Hessian equation, satisfying $\|\H y-r\|\le\epsilon\|r\|$.
\begin{enumerate}
\item (Formation) Compute the well-conditioned plus low-rank decomposition
(\ref{eq:hessExp2}).
\begin{enumerate}
\item Compute eigendecomposition $W=V\Lambda V^{T}$ and set $\tau=\lambda_{\min}(W)$.
\item Form the matrices $W_{0}$ and $U$ via (\ref{eq:Wsplit2}), and compute
the Cholesky factorization $ZZ^{T}=2W_{0}+UU^{T}$.
\end{enumerate}
\item (Factorization) Form the size-$nk$ Schur complement $\Sch=\tau I+(U\otimes Z)^{T}\A(\A^{T}\A)^{-1}\A^{T}(U\otimes Z)$
and compute its Cholesky factorization $\mathbf{L}\mathbf{L}^{T}=\mathbf{S}$.
\item (Solution) Use preconditioned conjugate gradients (PCG) to solve $\H y=r$
with $\tilde{\H}$ as preconditioner to $\epsilon$ relative residual.
Do at each PCG iteration:
\begin{enumerate}
\item Compute the matrix-vector product with $\H$ using the Kronecker identity
in (\ref{eq:hesseqn2}). 
\item Compute the matrix-vector product with $\tilde{\H}^{-1}$ using the
Sherman\textendash Morrison\textendash Woodbury in (\ref{eq:SMW2}),
implementing each $\Sch^{-1}=\mathbf{L}^{-T}\mathbf{L}^{-1}$.
\end{enumerate}
\end{enumerate}
\end{algor}

The main set-up cost is the factorization of the preconditioner (Step
2), which requires $nk$ matrix-vector products with $\A^{T}$, $\A$,
$(\A^{T}\A)^{-1}$, and $(U\otimes Z)^{T}$, and a single dense size-$nk$
Cholesky factorization. Under Assumption~\ref{ass:sparse}, this
requires
\[
(1/3)n^{3}k^{3}+O(n^{3}k^{2})\text{ flops and }\Theta(n^{2}k^{2})\text{ memory.}
\]
(Note that we have used $m\le n^{2}$.) The method converges to an
$\epsilon$-accurate solution in at most $\frac{1}{2}\kappa_{0}\log(2\kappa_{0}/\epsilon)$
PCG iterations, where $\kappa_{0}=\kappa(W_{0})$ as in Lemma~\ref{lem:Htilde},
and each iteration requires
\[
2n^{3}+n^{2}k^{2}+O(n^{2}k)\text{ flops.}
\]
The dominant $2n^{3}$ term arises from the matrix-vector product
$(W\otimes W)\vector X=\vector(WXW)$, as a part of the matrix-vector
product with $\H$. The $n^{2}k^{2}$ term arises from the application
of the Schur complement inverse $\Sch^{-1}=\mathbf{L}^{-T}\mathbf{L}^{-1}$.
Dropping the lower-order terms yields the following complexity estimate.
\begin{thm}
\label{thm:complex}Algorithm~\ref{alg:prec} uses $\Theta(n^{2}k^{2})$
memory and terminates after $\Theta(n^{3}k^{3}+n^{3}\log(1/\epsilon))$
flops.
\end{thm}
It is interesting to note that the complexity figure is not strongly
affected by the exact value of $m$. By comparison, explicitly forming
and factorizing the Hessian matrix $\H=\A^{T}(W\otimes W)\A$ under
Assumption~\ref{ass:sparse} requires
\begin{equation}
(1/3)m^{3}+O(n^{3}m+m^{2})\text{ flops and }\Theta(m^{2})\text{ memory.}\label{eq:complex_direct}
\end{equation}
Hence, our algorithm yields the biggest speed-up when the number of
constraints $m$ is large, and when the ratio $nk/m\ll1$. In particular,
it is up to a factor of $\sim n^{3}$ more efficient for problems
with number of constraints $m\sim n^{2}$. 

\subsection{Relation with prior work}

The CG (or PCG) solution of the interior-point Hessian equation is
an old idea that remains the standard approach for network-flow linear
programs~\cite[Ch.4]{mitchell1998interior}, and in general-purpose
solvers for nonlinear programming~\cite{byrd1999interior}; see also~\cite{benzi2005numerical}
and the references therein. The CG approach has not found widespread
use in SDP solvers, however, due to the considerable difficulty in
formulating an effective preconditioner. Existing preconditioners
had primarily been based on sparse matrix ideas, but these are not
applicable to the fully-dense Hessian equations arising from SDPs. 

Toh and Kojima~\cite{toh2002solving} were the first to develop highly
effective \emph{spectral} preconditioners based on the low-rank perturbed
view of the scaling matrix $W=W_{0}+UU^{T}$, but its use required
almost as much time and memory as a single iteration of the regular
interior-point method. Our preconditioner is similar in spirit, but
we make a number of modifications to improve efficiency. In particular,
our use of the Sherman\textendash Morrison\textendash Woodbury identity
allows us to prove a formal complexity bound that is strictly better
than the standard approach based on Cholesky factorization.

\section{\label{sec:num_stab}Improving Numerical Stability}

Unfortunately, the preconditioner in the previous section suffers
from numerical issues as the outer interior-point approaches the exact
solution. The culprit is the Sherman\textendash Morrison\textendash Woodbury
(SMW) formula, which is well-known to be numerically unstable when
the perturbed matrix is ill-conditioned; see e.g.~\cite{yip1986note}.

\subsection{Solving an augmented system}

Consider, for example, solving the preconditioner equation $\tilde{W}x=b$
from (\ref{eq:tIUUt}) at an interior-point step with duality gap
parameter $\mu$. The governing matrix $\tilde{W}=\tau I+UU^{T}$
becomes highly ill-conditioned as $\mu\to0^{+}$, with condition number
scaling $\kappa(\tilde{W})=\Theta(1/\mu)$. To avoid the SMW formula,
a standard implementation trick is to solve the symmetric indefinite
augmented problem
\begin{align}
\begin{bmatrix}\tau I & \sqrt{\tau}U\\
\sqrt{\tau}U^{T} & -\tau I
\end{bmatrix}\begin{bmatrix}x\\
y
\end{bmatrix} & =\begin{bmatrix}b\\
0
\end{bmatrix}.\label{eq:aug_ex}
\end{align}
Observe that performing Gaussian elimination (without pivoting) on
(\ref{eq:aug_ex}) results in \emph{identical} steps to a direct application
of the SMW formula (\ref{eq:wsmf}). However, the augmented system
is considerably better conditioned, with condition number $\sqrt{1+\|U\|^{2}/\tau}=\Theta(1/\sqrt{\mu})$.
This is a square-root factor better than $\tilde{W}$ itself, so we
would expect to lose half as many digits to round-off error as the
SMW formula by solving (\ref{eq:aug_ex}) using a stable method, like
LDL Cholesky factorization with numerical pivoting. In practice, numerical
pivoting usually results in some loss of efficiency. An acceptible
trade-off can generally be achieved by adjusting the ``threshold''
parameter for numerical pivots; see~e.g.~\cite{duff2004ma57}.

\subsection{An augmented preconditioner}

The augmented system approach cannot be directly applied to the preconditioner
$\tilde{\H}=\A^{T}\A+\U\U^{T}$, without considerably increasing the
cost of Algorithm~\ref{alg:prec}. This discrepency lies in the fact
that $\U$ is dense, containing $mnk$ nonzeros, but can be applied
in just $O(n^{2}k+m)$ flops using (\ref{eq:Umv}), as if it were
sparse. This special structure is lost when $\tilde{\H}$ is posed
in its augmented system form, and $\U$ is treated like any regular
dense matrix.

In the case that data matrix $\A$ is sparse, we may consider making
the following modification to $\tilde{\H}$:
\begin{equation}
\hat{\H}\triangleq\A^{T}\left(\tau^{2}I+(2\tau)UU^{T}\otimes I\right)\A,\label{eq:sparse_prec}
\end{equation}
which further approximates the dense matrix $2W_{0}+UU^{T}=ZZ^{T}$
using the scaled identity $ZZ^{T}\approx2\tau I$.
\begin{lem}
\label{lem:aug_prec}Let $\H$ and $\hat{\H}$ be defined in (\ref{eq:sparse_prec}),
and choose $\tau$ to satisfy $\lambda_{\min}(W_{0})\le\tau\le\lambda_{\max}(W_{0})$.
Then, $\lambda_{j}(\hat{\H}^{-1}\H)/\lambda_{n}(\hat{\H}^{-1}\H)\le\kappa^{2}(W_{0})$
for $j>k^{2}$.
\end{lem}
\begin{IEEEproof}
Define the $m\times nk$ matrix $F\triangleq\begin{bmatrix}I_{m} & \sqrt{2\tau}(U\otimes I_{n})^{T}\end{bmatrix}^{T}\A$,
the $nk\times k^{2}$ matrix $V=\begin{bmatrix}0_{nk} & I_{k}\otimes U\end{bmatrix}^{T}$,
and note that $\H_{p}\triangleq\hat{\H}^{-1/2}\H\hat{\H}^{-1/2}$
can be written $\H_{p}=Q^{T}\diag(W_{0}\otimes W_{0},\tau I\otimes W_{0})Q+\frac{\tau}{2}(Q^{T}V)(Q^{T}V)^{T}$
where $Q\triangleq F(F^{T}F)^{-1/2}$ is orthonormal. By the Cauchy
interlacing eigenvalues theorem, the first matrix has eigenvalues
that lie within the interval $\mathcal{I}\triangleq[\lambda_{\min}^{2}(W_{0}),\lambda_{\max}^{2}(W_{0})]$.
The second matrix is rank-$k^{2}$ and positive semidefinite, so can
perturb at most $k^{2}$ eigenvalues. We have $\lambda_{\min}^{2}(W_{0})\le\lambda_{j}(\H_{p})\le\lambda_{\max}^{2}(W_{0})$
for all $k^{2}<j\le m$, thereby yielding the desired result. 
\end{IEEEproof}
In view of Proposition~\ref{prop:PCG}, PCG with $\hat{\H}$ as preconditioner
converges to an $\epsilon$-accurate solution of the Hessian equation
$\H y=r$ in $k^{2}+O(\log\epsilon^{-1}\mu^{-1})$ iterations. The
figure is $O(\log(1/\epsilon))$ for all practical purposes, because
the ratio between $\epsilon$ and $\mu$ must be kept approximately
constant for the outer interior-point method to maintain its usual
convergence rate.

At each PCG iteration, the matrix-vector product with $\hat{\H}^{-1}$
may be implemented by solving the sparse augmented system 
\begin{equation}
\begin{bmatrix}\tau^{2}\A^{T}\A & \tau^{3/2}\A^{T}(U\otimes I)\\
\tau^{3/2}(U\otimes I)^{T}\A & -\tau^{2}/2I_{nk}
\end{bmatrix}\begin{bmatrix}x\\
y
\end{bmatrix}=\begin{bmatrix}b\\
0
\end{bmatrix}.\label{eq:aug_prec}
\end{equation}
The matrix condition number scales $\Theta(1/\sqrt{\mu})$, and some
bookkeeping shows that precomputing the LDL Cholesky without numerical
pivoting attains the same $O(n^{3}k^{3})$ factorization and $O(n^{2}k^{2})$
application costs as Algorithm~\ref{alg:prec} in Theorem~\ref{thm:complex}.
If the matrix sparsity pattern of (\ref{eq:aug_prec}) is structured
in a nice way, then it is often possible for a sparse factorization
of (\ref{eq:aug_prec}) to be computed at even further reduced costs.
Indeed, the matrix contains just $\sim mk+nk$ nonzeros, so the cost
of sparse Cholesky factorization can be as low as $\sim n^{2}k^{2}$,
or even as $\sim nk$.

\section{\label{sec:Numerical-Results}Numerical Results}

\begin{figure}
\hfill{}\includegraphics[width=0.8\columnwidth]{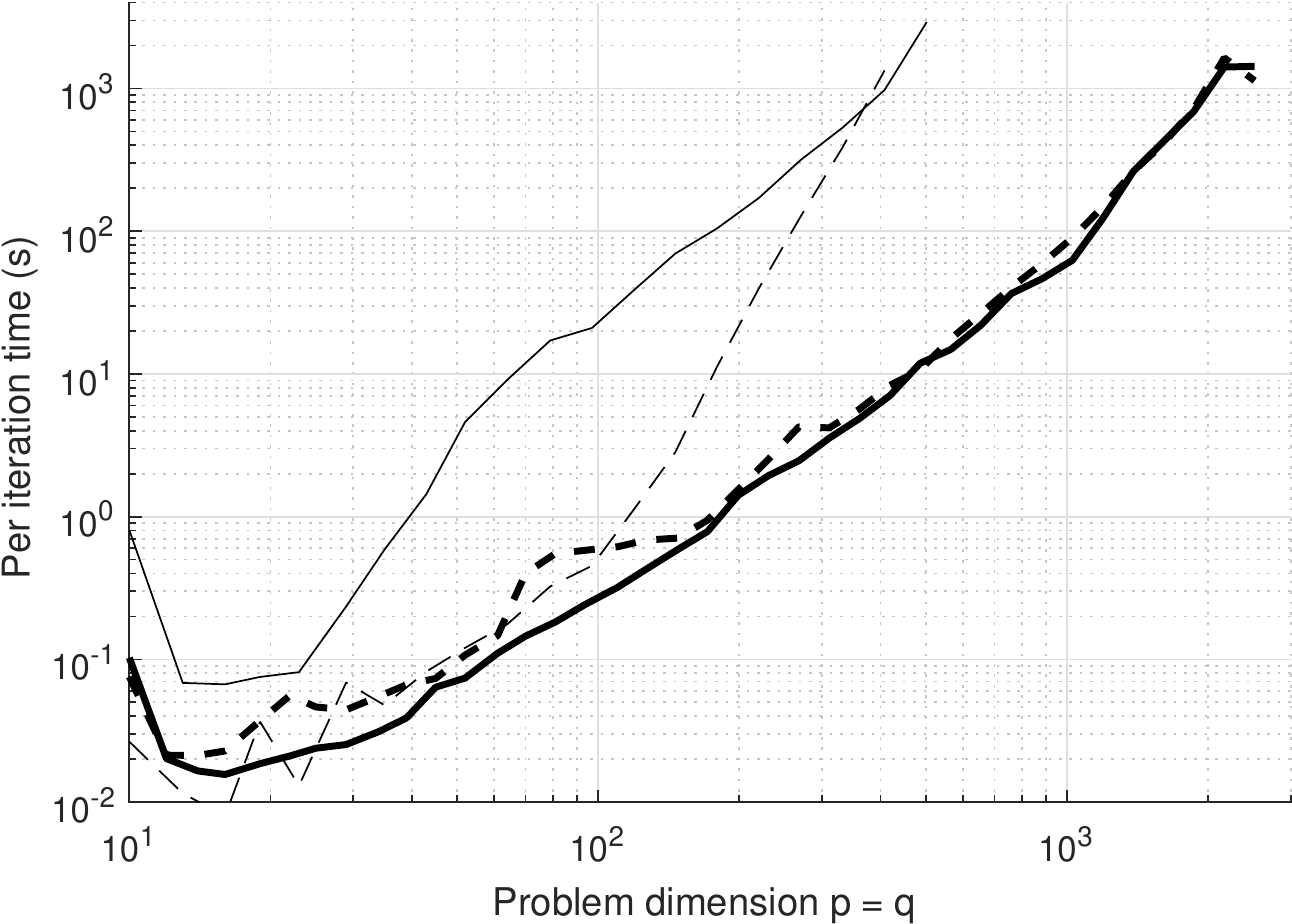}\hfill{}

\caption{\label{fig:per_iter}Per interior-point iteration time in seconds
for modified SeDuMi (thick lines) and regular SeDuMi (thin lines)
for matrix completion SDPs with $p=q$, rank $k=1$ and: (solid) $m=50p=25n$
constraints; (dashed) $m=0.1pq=0.025n^{2}$ constraints.}
\end{figure}
\begin{figure}
\hfill{}\includegraphics[width=0.8\columnwidth]{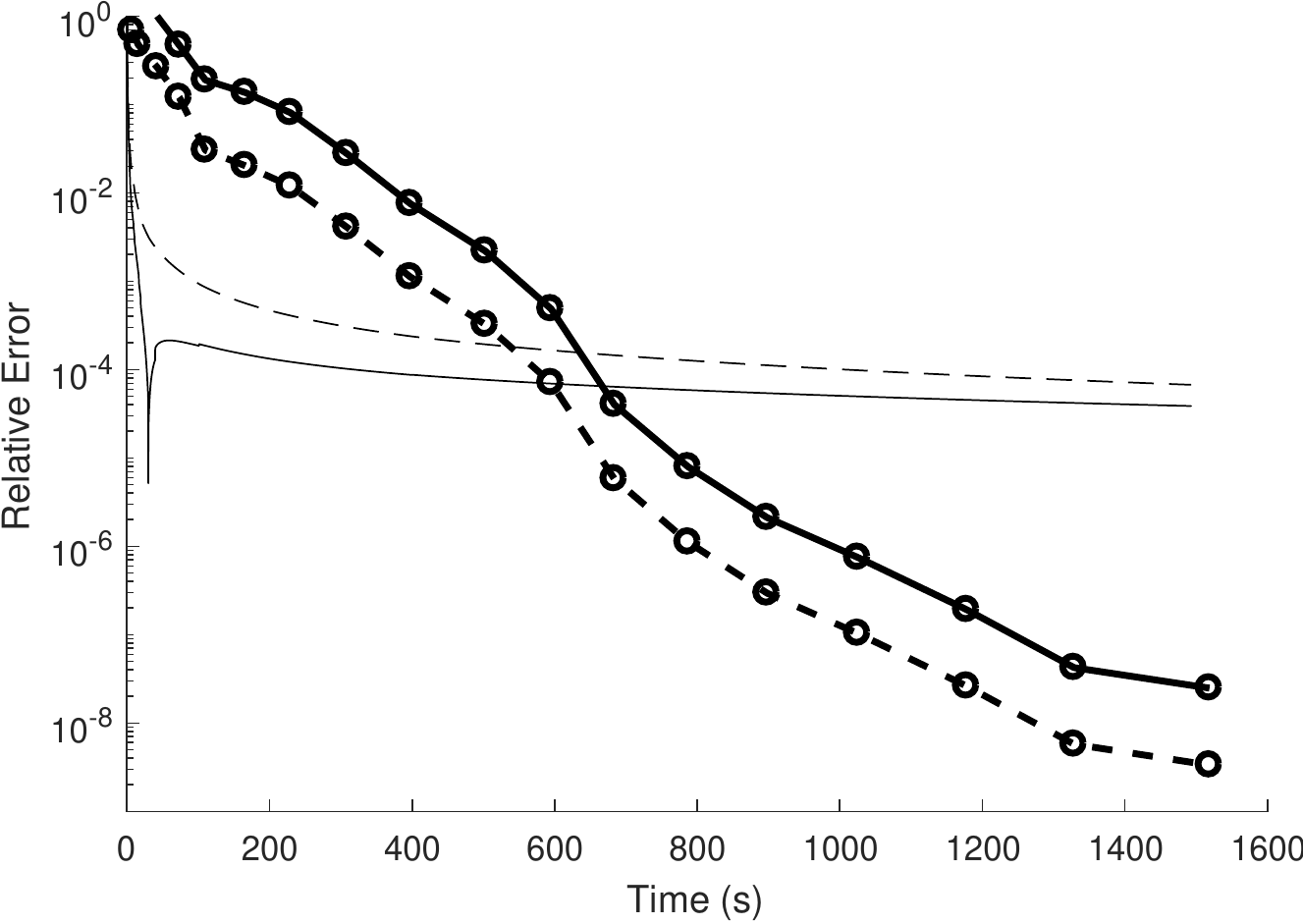}\hfill{}

\caption{\label{fig:err_rate}Progress of modified SeDuMi (thick lines) and
singular value thresholding (thin lines) for a matrix completion problem
with $p=q=500$, rank $k=4$, and $m=20,000$ constraints: (solid)
objective error $\mathrm{abs}(\|Z\|_{*}-\|M\|_{*})/\|M\|_{*}$; (dashed)
relative residual $\left(\sum_{i,j\in\Omega}(Z_{i,j}-M_{i,j})^{2}/\sum_{i,j\in\Omega}M_{i,j}^{2}\right)^{1/2}$.
Each dot represents a single interior-point iteration. }
\end{figure}
We implemented the preconditioner in Section~\ref{sec:num_stab}
in MATLAB, and embedded it within SeDuMi version 1.3~\cite{sturm1999sedumi};
the resulting solver is publicly available at
\[
\text{{\tt http://alum.mit.edu/www/ryz}}
\]
SeDuMi is an NT-scaled strictly feasible path-following interior-point
method, so we expect all of our theoretical results to hold. In fact,
the original SeDuMi code already incorporates PCG in its solution
of the Hessian equation, but the preconditioner is a numerically stabilized
Cholesky factorization of the actual Hessian matrix $\H=\A^{T}\D\A$.
Therefore, our only substantial modification is to replace this near-exact
preconditioner with the spectral approximation $\hat{\H}$, implemented
using the augmented system representation in (\ref{eq:aug_prec}).
The LDL Cholesky factorization is computed using the \texttt{ldl}
command in MATLAB, which calls the MA57 routine by Duff~\cite{duff2004ma57}.

For a general SDP, the exact value of $k=\rank X^{\star}$ is unknown
until after the problem has already been solved. In theory, our guarantees
will continue to hold by setting the rank parameter to any upper-bound
$k_{\max}=O(k)$, but in practice, the algorithm will run considerably
faster using a less pessimistic value. Our implementation uses the
spectrum of the scaling matrix $W$ to dynamically estimate a reasonable
approximation $\tilde{k}\approx k$. More specifically, given an upper-bound
$k_{\max}\ge k$ and an eigenvalue ratio $\eta$, we set $\tilde{k}$
as:
\[
\tilde{k}=\max\{i\in\{0,1,\ldots,k_{\max}\}:\lambda_{i}(W)\ge\eta\lambda_{i+1}(W)\}.
\]
The heuristic is inspired by Proposition~\ref{prop:Weig}: as the
interior-point method progresses and the duality gap parameter $\mu\to0^{+}$,
the true value of $k$ causes the ratio $\lambda_{k}(W)/\lambda_{k+1}(W)\in\Theta(1/\mu)$
to tend to infinity. In other words, $\tilde{k}$ is guaranteed to
converge to the true $k$ as the interior-point method progresses
towards the solution. 

\subsection{Test problem: Matrix completion}

The matrix completion problem seeks to recover a low-rank size-$p\times q$
rectangular matrix $M$, by observing an incomplete subset of entries
$M_{i,j}$ at $\{i,j\}\in\Omega$ and solving the convex optimization
program
\begin{equation}
Z^{\star}=\text{ minimize }\|Z\|_{*}\text{ s.t. }Z_{i,j}=M_{i,j}\;\forall\{i,j\}\in\Omega,\label{eq:matcomp}
\end{equation}
where the nuclear norm $\|Z\|_{*}=\tr(Z^{T}Z)^{1/2}$ is the sum of
the singular values. Note that (\ref{eq:matcomp}) is a size $n=p+q$
semidefinite program over $m=|\Omega|$ constraints
\begin{align}
\text{ minimize } & I\bullet X\label{eq:matcomp_sdp}\\
\text{subject to } & \frac{1}{2}\begin{bmatrix}0_{p} & E_{i,j}^{T}\\
E_{i,j} & 0_{q}
\end{bmatrix}\bullet X=M_{i,j}\;\forall\{i,j\}\in\Omega\nonumber \\
 & \begin{bmatrix}U & Z^{T}\\
Z & V
\end{bmatrix}=X\succeq0,\nonumber 
\end{align}
where $E_{i,j}$ is an $p\times q$ matrix containing a single ``1''
at its $\{i,j\}$-th element. It is a famous result by Candes and
Recht~\cite{candes2012exact}, later improved by Candes and Tao~\cite{candes2010power}
that, when $M$ is low-rank and incoherent, and the number of samples
satisfy $m\ge Cn(\log n)^{2}$ with some constant $C$, then all $pq$
elements of $M$ are exactly recovered by solving (\ref{eq:matcomp}).
In other words, the solution to (\ref{eq:matcomp}) is precisely $Z^{\star}=M$. 

Matrix completion makes an ideal test problem for the PCG procedure
described in this paper, for the following reasons: 
\begin{enumerate}
\item The solution rank $k=\rank X^{\star}$ is easily adjustable by controlling
the rank of the original matrix $M$;
\item The SDP order $n=p+q$ and the number of constraints $m$ are easily
adjustable by controlling the size of the original matrix $M$ and
by modifying the number of observations $|\Omega|$;
\item The data matrix $\A=[\vector A_{1},\ldots,\vector A_{m}]$ is highly
sparse, and always satisfies Assumption~\ref{ass:sparse} by construction. 
\end{enumerate}
In this section, we consider random instances of (\ref{eq:matcomp_sdp}).
More specifically, we select $\Omega\subseteq\{1,\ldots,p\}\times\{1,\ldots,q\}$
uniformly at random from all subsets with cardinality $m$, and set
$M=G_{1}G_{2}^{T}$, where $G_{1}\in\R^{p\times k}$ and $G_{2}\in\R^{q\times k}$
are selected i.i.d. from the standard Gaussian.

\subsection{Comparison with standard SeDuMi}

The matrix completion SDP (\ref{eq:matcomp_sdp}) is a challenging
test problem for all standard interior-point solvers. The bottleneck
is factoring the $m\times m$ fully-dense Hessian matrix $\H$, for
worst-case complexities of $O(n^{6})$ time and $O(n^{4})$ memory.
Chordal decomposition cannot be used to reduce these complexity figures,
because the underlying graph does not have a bounded treewidth; see~\cite{gao2012treewidth}. 

By comparison, our modified SeDuMi gains considerable efficiency by
avoiding an explicit treatment of the Hessian matrix $\H$. In all
of our numerical trials, the augmented system (\ref{eq:aug_prec})
associated with the preconditioner $\hat{\H}$ is highly sparse, and
the algorithm's bottleneck is the matrix-vector product $(W\otimes W)\vector X=\vector(WXW)$,
as a part of the matrix-vector product with $\H$. These are realized
as matrix-matrix products and evaluated using BLAS routines, so our
MATLAB implementation should have a comparable level of performance
to a hand-coded C/C++ implementation.

Figure~\ref{fig:per_iter} compares the per-iteration cost of our
modified SeDuMi and the standard implementation, on a modest workstation
with 16 GB of RAM and an Intel Xeon E5-2609 v4 CPU with eight 1.70
GHz cores. Two sets of problems were considered: one set with $m=25n$
and another with $m=0.025n^{2}$. As shown, the time complexity of
standard SeDuMi is highly dependent upon the number of constraints
$m$, but this dependency is essentially eliminated in the modified
version. Standard SeDuMi was able to solve problems with $p+q=n\approx800$
before running of memory. By contrast, our modified SeDuMi was able
to solve matrix completion problems as large as $n=5024$ and $m=6.31\times10^{5}$,
in around 8 hours. Simply storing the associated Hessian matrix would
have required 1,600 GB of memory, which is a hundred times what was
available. In all of these trials, PCG converges to an iterate of
sufficient accuracy in 15-25 iterations (except when stagnation occurs
due to numerical issues). 

\subsection{Comparison with singular value thresholding}

Our modified SeDuMi is a true second-order method, because it converges
at a linear rate, requiring $O(\log(1/\epsilon))$ iterations to produce
an $\epsilon$-accurate solution. To make this distinction clear,
we compare our modified SeDuMi method with the singular value thresholding
(SVT) algorithm, a popular and widely-used first-order method for
matrix-completion problems~\cite{cai2010singular}. The SVT algorithm
implicitly represents $Z$ in its low-rank factored form, and computes
singular values using the Lanczos iteration; its per-iteration complexity
is as low as $\sim mk$ time and $\sim nk+m$ memory. However, the
method converges sublinearly in the worst-case, requiring $O(1/\epsilon)$
iterations to produce an $\epsilon$-accurate solution. 

Figure~\ref{fig:err_rate} shows the progress of our modified SeDuMi
and SVT over a 25 minute period, for a random matrix completion problem
with $p=q=500$, rank $k=4$, and $m=2\times10^{4}$ observations.
After 20,000 iterations, SVT outputs an estimation of $M$ with relative
error of $\approx10^{-4}$. Indeed, SVT was able to compute an iterate
of nearly this accuracy in just 2 minutes, but its sublinear convergence
rate produces diminishing returns for the additional computation time.
By contrast, modified SeDuMi converges linearly, gaining one decimal
digit of accuracy every 3 minutes. After 18 outer interior-point iterations
and 4233 inner PCG iterations, the method outputs an estimation of
$M$ with relative error of $\approx10^{-8}$.

\section{Conclusion}

This paper describes a preconditioner that allows preconditioned conjugate
gradients (PCG) to converge to a solution of the interior-point Hessian
equation in a few tens of iterations, independent of the ill-conditioning
of the Hessian matrix. The preconditioner can be factored in $\Theta(n^{3}k^{3})$
time and $\Theta(n^{2}k^{2})$ memory, and the cost of the subsequent
PCG iterations is dominated by matrix-vector products with the Hessian
matrix. We embed the preconditioner within SeDuMi, and use it to solve
large-and-sparse, low-rank, matrix completion SDPs to 8-10 decimal
digits of accuracy. The largest problem we considered had $n=5024$
and $m=6.31\times10^{5}$, and was solved in less than 8 hours on
a modest workstation with 16 GB of memory.

\bibliographystyle{IEEEtran}
\bibliography{lowrank}

\end{document}